\documentclass[a4paper,dvipdfmx]{article}
\usepackage{a4wide}
\usepackage[dvipdfmx]{graphicx}
\usepackage[dvipdfmx]{color}
\usepackage{color}
\usepackage{bm}
\usepackage{amsmath}
\usepackage{amssymb}
\usepackage{amsfonts}
\usepackage{amsmath}
\usepackage{ascmac}
\usepackage[english]{babel}
\usepackage{mathrsfs}
\title{On the automorphism of Barns Wall Lattice $\Lambda_{BW_{16}}$ and rank 4 tensor of quaternions}
\author{Misaki Ohta,University of the Ryukyus \thanks{Department of physics, Email: e193225@eve.u-ryukyu.ac.jp or apple.designed@icloud.com}}
\begin{document}
\maketitle

\begin{abstract}
In a previous paper, I found that the Weyl group $W(F_4)$ and Barns-Wall Lattice $BW_{16}$ can be constructed using the rank $2$ tensor of the quaternion. In the present paper, I describe how I were able to construct an algebra, which is the subalgebra of the direct product of Hurwitz Quaternionic integers $\mathscr{H}^4$,  isomorphic to the automorphism $\text{Aut}(BW_{16})$ order $2^{21} \cdot 3^5 \cdot 5^2 \cdot 7$ of Barns Wall Lattice $BW_{16}$ by functionally extending the rank of the tensor product of quaternions to $4$.
\end{abstract}

\section{The rank 1 tensor of the Quaternions and alternating group $2 \cdot A_4$}
Quaternions $\mathbb{H}$ is defined by
$$
\mathbb{H}=\left\{z=x_0+x_1 i+x_2 j+x_3 k \mid x_i \in \mathbb{R}\right\}
$$
where
$$
\begin{aligned}
&i^{2}=j^{2}=k^{2}=i j k=-1 \quad (-e) \\
&i j=-j i=k, \quad j k=-k j=i, \quad k i=-i k=j
\end{aligned}
$$
And there is the subalgebra called Hurwitz Quaternionic integers
$$
\mathscr{H}=\{a+i b+j c+k d| a, b, c, d \in (\mathbb{Z} \oplus \mathbb{Z}+1 / 2)\}
$$
According to J.H.Conway and Neal.J.A.Sloane[2], the multiplicative group generated below, whose norm is 1, is isomorphic to $2 \cdot A_4$ order $24$.
$$
<i, j, k, \omega>
$$
where $\omega=\frac{1+i+j+k}{2}$.\\\\
In addition, all the elements of this group are transferred to the root system $F_4$ lattice by the following inverse mapping.
$$
\begin{array}{llllll}
\phi: &(\mathbb{Z}\oplus \mathbb{Z}+1/2)^4& \rightarrow& <i,j,k,\omega>&\\
&(c_0,c_1,c_2,c_3)& \mapsto &  c_0 + c_1 i + c_2 j + c_3 k
\end{array}
$$
for example
$$
\begin{array}{lllll}
&\phi^{-1}: &i &\mapsto &(0,1,0,0)\\
&\phi^{-1}: &j &\mapsto &(0,0,1,0)\\
&\phi^{-1}: &k &\mapsto &(0,0,0,1)\\
&\phi^{-1}: &\omega &\mapsto &(1/2,1/2,1/2,1/2)\\
\end{array}
$$
The $F_4$ Lattice is as follows.
$$
\Lambda_{F_4}=\text { RowSpan }_{\mathbf{Z}}\frac{1}{2}
\left(\begin{array}{cccc}
1 & 1 & 1 & 1 \\
& 2 & & \\
& & 2 & \\
& & & 2
\end{array}\right)
$$
The order of the automorphism of the lattice $W(F_4)$ is $1152$ and the kissing number of the lattce is $24$.
\section{The rank 2 tensor of the Quaternions and Weyl Group $W(F_4)$}
For $a_i,b_j \in \{\pm e, \pm i , \pm j, \pm k \}$, it is clear that the tensor product satisfies
$$
\left(a_{1} \otimes b_{1}\right)\left(a_{2} \otimes b_{2}\right)=a_{1} a_{2} \otimes b_{1} b_{2}
$$
The exceptional Weyl group $G = W(F_4)$, whose order is $1152 = 24 \times 24 \times 2$, is generated by
$$
\begin{aligned}
e_1 &= \frac{1}{2}(e \otimes e+i \otimes j+j \otimes k-  k \otimes i)\\
e_2 &= \frac{1}{2}(e \otimes e-i \otimes i-j \otimes j-  k \otimes k)\\
e_3 &= \frac{1}{2}(e \otimes e-i \otimes i-j \otimes k+  k \otimes j)\\
e_4 &= \frac{1}{2}(e \otimes e-i \otimes j-j \otimes i+  k \otimes k)\\
\end{aligned}
$$
where
$$
\begin{aligned}
&i^{2}=j^{2}=k^{2}=i j k=-1 \quad (-e) \\
&i j=-j i=k, \quad j k=-k j=i, \quad k i=-i k=j
\end{aligned}
$$
This is a subalgebra whose norm is 1 in the direct product $\mathscr{H}^2$ of the algebra called Hurwitz Quaternionic integers defined below.
$$
\mathscr{H}=\{a+i b+j c+k d| a, b, c, d \in (\mathbb{Z} \oplus \mathbb{Z}+1 / 2)\}
$$
At the same time,  $G$ forms a Barnes-Wall lattice $\Lambda$ defined by the following in the same way as for the $F_4$ lattice(See the previous paper[1]).
$$
\Lambda = \operatorname{RowSpan}_{\mathbf{Z}}\frac{1}{4}\left(\begin{array}{cccccccccccccccc}
1 & 1 & 1 & 1 & 1 & 1 & 1 & 1 & 1 & 1 & 1 & 1 & 1 & 1 & 1 & 1 \\
 & 2 &  &  &  &  &  & 2 &  &  &  & 2 &  & 2 &  &  \\
 &  & 2 &  &  &  &  & 2 &  &  &  & 2 &  &  & 2 &  \\
 &  &  & 2 &  &  &  & 2 &  &  &  & 2 &  &  &  & 2 \\
 &  &  &  & 2 &  &  & 2 &  &  &  &  &  & 2 & 2 &  \\
 &  &  &  &  & 2 &  & 2 &  &  &  &  &  & 2 &  & 2 \\
 &  &  &  &  &  & 2 & 2 &  &  &  &  &  &  & 2 & 2 \\
 &  &  &  &  &  &  & 4 &  &  &  &  &  &  &  &  \\
 &  &  &  &  &  &  &  & 2 &  &  & 2 &  & 2 & 2 &  \\
 &  &  &  &  &  &  &  &  & 2 &  & 2 &  & 2 &  & 2 \\
 &  &  &  &  &  &  &  &  &  & 2 & 2 &  &  & 2 & 2 \\
 &  &  &  &  &  &  &  &  &  &  & 4 &  &  &  &  \\
 &  &  &  &  &  &  &  &  &  &  &  & 2 & 2 & 2 & 2 \\
 &  &  &  &  &  &  &  &  &  &  &  &  & 4 &  &  \\
 &  &  &  &  &  &  &  &  &  &  &  &  &  & 4 &  \\
 &  &  &  &  &  &  &  &  &  &  &  &  &  &  & 4
\end{array}\right)
$$
The order of the automorphism of this lattice is $89181388800=2^{21} \cdot 3^5 \cdot 5^2 \cdot 7$ and its stracture is $\text{Aut}(\Lambda)\cong 2^{1+8} \mathrm{PSO}_8{ }^{+}\left(\mathrm{F}_2\right)$, the kissing number is $4320$ as well.

\section{The rank 4 tensor of the Quaternions and the group $\operatorname{Aut}\left(\Lambda_{\text {Barns-Wall }}\right)$}
I now wish to consider the frank 4 tensor of the quaternion, which is a subalgebra of $\mathscr{H}^4$.

\begin{itembox}{Fact.1}
The group generated by the following elements is isomorphic to $\operatorname{Aut}\left(\Lambda_{\text {Barns-Wall }}\right)$ order $89181388800=2^{21} \cdot 3^5 \cdot 5^2 \cdot 7$.
$$
\begin{aligned}
x_1=&\frac{1}{2}[(e \otimes e \otimes e \otimes e)+(i \otimes j \otimes k \otimes i)+(j \otimes k \otimes i \otimes k)+(k \otimes i \otimes j \otimes j)]\\
x_2=&\frac{1}{2}[(e \otimes e \otimes e \otimes e)+(i \otimes i \otimes k \otimes i)+(j \otimes k \otimes i \otimes k)-(k \otimes j \otimes j \otimes j)]\\
x_3=&\frac{1}{2}[(e \otimes e \otimes e \otimes e)+(j \otimes j \otimes i \otimes i)+(i \otimes k \otimes k \otimes k)+(k \otimes i \otimes j \otimes j)]\\
x_4=&\frac{1}{2}[(e \otimes e \otimes e \otimes e)+(k \otimes j \otimes i \otimes k)+(i \otimes k \otimes k \otimes j)-(j \otimes i \otimes j \otimes i)]\\
x_5=&\frac{1}{2}[(e \otimes e \otimes e \otimes e)+(k \otimes j \otimes k \otimes k)-(i \otimes k \otimes i \otimes j)-(j \otimes i \otimes j \otimes i)]\\
x_6=&\frac{1}{2}[(e \otimes e \otimes e \otimes e)+(i \otimes i \otimes e \otimes e)+(j \otimes k \otimes e \otimes e)+(k \otimes j \otimes e \otimes e)]\\
x_7=&\frac{1}{2}[(e \otimes e \otimes e \otimes e)+(j \otimes j \otimes k \otimes e)+(i \otimes k \otimes i \otimes e
)+(k \otimes i \otimes j \otimes e)]
\end{aligned}
$$
\end{itembox}
You can easily verify the above fact Using MAGMA and the following representation. The representation actually satisfies the following.
$$
\begin{aligned}
\rho(\left(y_{1} \otimes e \otimes e \otimes e\right)) \cdot \rho(\left(e \otimes y_{2} \otimes e \otimes e\right))\\
=\rho(\left(e \otimes y_2 \otimes e \otimes e\right)) \cdot \rho(\left(y_1 \otimes e \otimes e \otimes e\right))
\end{aligned}
$$

$$
\begin{aligned}
\rho(\left(y_{1} \otimes e \otimes e \otimes e\right)) \cdot \rho(\left(e \otimes e \otimes y_{3} \otimes e\right))\\
=\rho((\left(e \otimes e \otimes y_{3} \otimes e\right)) \cdot \rho(\left(y_1 \otimes e \otimes e \otimes e\right))
\end{aligned}
$$

$$
\vdots
$$

$$
\begin{aligned}
\rho(\left(e \otimes e \otimes y_{3} \otimes e\right)) \cdot \rho(\left(e \otimes e \otimes e \otimes y_{4}\right))\\
=\rho(\left(e \otimes e \otimes e \otimes y_{4} \right)) \cdot \rho(\left(e \otimes e \otimes y_{3} \otimes e\right))
\end{aligned}
$$
and
$$
\rho(\left(y_{1} \otimes e \otimes e \otimes e\right))^2= \rho(\left(e \otimes y_{2}\otimes e \otimes e\right))^2=\rho(\left(e \otimes e\otimes y_3 \otimes e\right))^2=\rho(\left(e \otimes e\otimes e \otimes y_4\right))^2 =- \rho((e\otimes e \otimes e \otimes e))
$$
where $y_i \in \{i,j,k=ij\}$\\
******************************************************************************************
$$
\rho((i \otimes e  \otimes e  \otimes e))=
\left(\begin{array}{rrrr|rrrr|rrrr|rrrr}
 & 1 &  &  &  &  &  &  &  &  &  &  &  &  &  &  \\
-1 &  &  &  &  &  &  &  &  &  &  &  &  &  &  &  \\
 &  &  & 1 &  &  &  &  &  &  &  &  &  &  &  &  \\
 &  & -1 &  &  &  &  &  &  &  &  &  &  &  &  &  \\
\hline
  &  &  &  &  & 1 &  &  &  &  &  &  &  &  &  &  \\
 &  &  &  & -1 &  &  &  &  &  &  &  &  &  &  &  \\
 &  &  &  &  &  &  & 1 &  &  &  &  &  &  &  &  \\
 &  &  &  &  &  & -1 &  &  &  &  &  &  &  &  &  \\
\hline
  &  &  &  &  &  &  &  &  & 1 &  &  &  &  &  &  \\
 &  &  &  &  &  &  &  & -1 &  &  &  &  &  &  &  \\
 &  &  &  &  &  &  &  &  &  &  & 1 &  &  &  &  \\
 &  &  &  &  &  &  &  &  &  & -1 &  &  &  &  &  \\
\hline
  &  &  &  &  &  &  &  &  &  &  &  &  & 1 &  &  \\
 &  &  &  &  &  &  &  &  &  &  &  & -1 &  &  &  \\
 &  &  &  &  &  &  &  &  &  &  &  &  &  &  & 1 \\
 &  &  &  &  &  &  &  &  &  &  &  &  &  & -1 &
\end{array}\right)
$$
$$
\rho((j \otimes e  \otimes e  \otimes e))=
\left(\begin{array}{rrrr|rrrr|rrrr|rrrr}
 &  & 1 &  &  &  &  &  &  &  &  &  &  &  &  &  \\
 &  &  & -1 &  &  &  &  &  &  &  &  &  &  &  &  \\
-1 &  &  &  &  &  &  &  &  &  &  &  &  &  &  &  \\
 & 1 &  &  &  &  &  &  &  &  &  &  &  &  &  &  \\
\hline
  &  &  &  &  &  & 1 &  &  &  &  &  &  &  &  &  \\
 &  &  &  &  &  &  & -1 &  &  &  &  &  &  &  &  \\
 &  &  &  & -1 &  &  &  &  &  &  &  &  &  &  &  \\
 &  &  &  &  & 1 &  &  &  &  &  &  &  &  &  &  \\
\hline
  &  &  &  &  &  &  &  &  &  & 1 &  &  &  &  &  \\
 &  &  &  &  &  &  &  &  &  &  & -1 &  &  &  &  \\
 &  &  &  &  &  &  &  & -1 &  &  &  &  &  &  &  \\
 &  &  &  &  &  &  &  &  & 1 &  &  &  &  &  &  \\
\hline
  &  &  &  &  &  &  &  &  &  &  &  &  &  & 1 &  \\
 &  &  &  &  &  &  &  &  &  &  &  &  &  &  & -1 \\
 &  &  &  &  &  &  &  &  &  &  &  & -1 &  &  &  \\
 &  &  &  &  &  &  &  &  &  &  &  &  & 1 &  &
\end{array}\right)
$$
$$
\rho((e \otimes i  \otimes e  \otimes e))=
\left(\begin{array}{rrrr|rrrr|rrrr|rrrr}
 &  &  &  & 1 &  &  &  &  &  &  &  &  &  &  &  \\
 &  &  &  &  & 1 &  &  &  &  &  &  &  &  &  &  \\
 &  &  &  &  &  & 1 &  &  &  &  &  &  &  &  &  \\
 &  &  &  &  &  &  & 1 &  &  &  &  &  &  &  &  \\
\hline
 -1 &  &  &  &  &  &  &  &  &  &  &  &  &  &  &  \\
 & -1 &  &  &  &  &  &  &  &  &  &  &  &  &  &  \\
 &  & -1 &  &  &  &  &  &  &  &  &  &  &  &  &  \\
 &  &  & -1 &  &  &  &  &  &  &  &  &  &  &  &  \\
\hline
  &  &  &  &  &  &  &  &  &  &  &  & 1 &  &  &  \\
 &  &  &  &  &  &  &  &  &  &  &  &  & 1 &  &  \\
 &  &  &  &  &  &  &  &  &  &  &  &  &  & 1 &  \\
 &  &  &  &  &  &  &  &  &  &  &  &  &  &  & 1 \\
\hline
  &  &  &  &  &  &  &  & -1 &  &  &  &  &  &  &  \\
 &  &  &  &  &  &  &  &  & -1 &  &  &  &  &  &  \\
 &  &  &  &  &  &  &  &  &  & -1 &  &  &  &  &  \\
 &  &  &  &  &  &  &  &  &  &  & -1 &  &  &  &
\end{array}\right)
$$
$$
\rho((e \otimes j  \otimes e  \otimes e))=
\left(\begin{array}{rrrr|rrrr|rrrr|rrrr}
 &  &  &  &  &  &  &  & 1 &  &  &  &  &  &  &  \\
 &  &  &  &  &  &  &  &  & 1 &  &  &  &  &  &  \\
 &  &  &  &  &  &  &  &  &  & 1 &  &  &  &  &  \\
 &  &  &  &  &  &  &  &  &  &  & 1 &  &  &  &  \\
\hline
  &  &  &  &  &  &  &  &  &  &  &  & -1 &  &  &  \\
 &  &  &  &  &  &  &  &  &  &  &  &  & -1 &  &  \\
 &  &  &  &  &  &  &  &  &  &  &  &  &  & -1 &  \\
 &  &  &  &  &  &  &  &  &  &  &  &  &  &  & -1 \\
\hline
 -1 &  &  &  &  &  &  &  &  &  &  &  &  &  &  &  \\
 & -1 &  &  &  &  &  &  &  &  &  &  &  &  &  &  \\
 &  & -1 &  &  &  &  &  &  &  &  &  &  &  &  &  \\
 &  &  & -1 &  &  &  &  &  &  &  &  &  &  &  &  \\
\hline
  &  &  &  & 1 &  &  &  &  &  &  &  &  &  &  &  \\
 &  &  &  &  & 1 &  &  &  &  &  &  &  &  &  &  \\
 &  &  &  &  &  & 1 &  &  &  &  &  &  &  &  &  \\
 &  &  &  &  &  &  & 1 &  &  &  &  &  &  &  &
\end{array}\right)
$$

$$
\rho((e \otimes e  \otimes i  \otimes e))=
\left(\begin{array}{rrrr|rrrr|rrrr|rrrr}
 & 1 &  &  &  &  &  &  &  &  &  &  &  &  &  &  \\
-1 &  &  &  &  &  &  &  &  &  &  &  &  &  &  &  \\
 &  &  & -1 &  &  &  &  &  &  &  &  &  &  &  &  \\
 &  & 1 &  &  &  &  &  &  &  &  &  &  &  &  &  \\
\hline
  &  &  &  &  & 1 &  &  &  &  &  &  &  &  &  &  \\
 &  &  &  & -1 &  &  &  &  &  &  &  &  &  &  &  \\
 &  &  &  &  &  &  & -1 &  &  &  &  &  &  &  &  \\
 &  &  &  &  &  & 1 &  &  &  &  &  &  &  &  &  \\
\hline
  &  &  &  &  &  &  &  &  & 1 &  &  &  &  &  &  \\
 &  &  &  &  &  &  &  & -1 &  &  &  &  &  &  &  \\
 &  &  &  &  &  &  &  &  &  &  & -1 &  &  &  &  \\
 &  &  &  &  &  &  &  &  &  & 1 &  &  &  &  &  \\
\hline
  &  &  &  &  &  &  &  &  &  &  &  &  & 1 &  &  \\
 &  &  &  &  &  &  &  &  &  &  &  & -1 &  &  &  \\
 &  &  &  &  &  &  &  &  &  &  &  &  &  &  & -1 \\
 &  &  &  &  &  &  &  &  &  &  &  &  &  & 1 &
\end{array}\right)
$$

$$
\rho((e \otimes e  \otimes j  \otimes e))=
\left(\begin{array}{rrrr|rrrr|rrrr|rrrr}
 &  & 1 &  &  &  &  &  &  &  &  &  &  &  &  &  \\
 &  &  & 1 &  &  &  &  &  &  &  &  &  &  &  &  \\
-1 &  &  &  &  &  &  &  &  &  &  &  &  &  &  &  \\
 & -1 &  &  &  &  &  &  &  &  &  &  &  &  &  &  \\
\hline
  &  &  &  &  &  & 1 &  &  &  &  &  &  &  &  &  \\
 &  &  &  &  &  &  & 1 &  &  &  &  &  &  &  &  \\
 &  &  &  & -1 &  &  &  &  &  &  &  &  &  &  &  \\
 &  &  &  &  & -1 &  &  &  &  &  &  &  &  &  &  \\
\hline
  &  &  &  &  &  &  &  &  &  & 1 &  &  &  &  &  \\
 &  &  &  &  &  &  &  &  &  &  & 1 &  &  &  &  \\
 &  &  &  &  &  &  &  & -1 &  &  &  &  &  &  &  \\
 &  &  &  &  &  &  &  &  & -1 &  &  &  &  &  &  \\
\hline
  &  &  &  &  &  &  &  &  &  &  &  &  &  & 1 &  \\
 &  &  &  &  &  &  &  &  &  &  &  &  &  &  & 1 \\
 &  &  &  &  &  &  &  &  &  &  &  & -1 &  &  &  \\
 &  &  &  &  &  &  &  &  &  &  &  &  & -1 &  &
\end{array}\right)
$$

$$
\rho((e \otimes e  \otimes e  \otimes i)) =
\left(\begin{array}{rrrr|rrrr|rrrr|rrrr}
 &  &  &  & 1 &  &  &  &  &  &  &  &  &  &  &  \\
 &  &  &  &  & 1 &  &  &  &  &  &  &  &  &  &  \\
 &  &  &  &  &  & 1 &  &  &  &  &  &  &  &  &  \\
 &  &  &  &  &  &  & 1 &  &  &  &  &  &  &  &  \\
\hline
 -1 &  &  &  &  &  &  &  &  &  &  &  &  &  &  &  \\
 & -1 &  &  &  &  &  &  &  &  &  &  &  &  &  &  \\
 &  & -1 &  &  &  &  &  &  &  &  &  &  &  &  &  \\
 &  &  & -1 &  &  &  &  &  &  &  &  &  &  &  &  \\
\hline
  &  &  &  &  &  &  &  &  &  &  &  & -1 &  &  &  \\
 &  &  &  &  &  &  &  &  &  &  &  &  & -1 &  &  \\
 &  &  &  &  &  &  &  &  &  &  &  &  &  & -1 &  \\
 &  &  &  &  &  &  &  &  &  &  &  &  &  &  & -1 \\
\hline
  &  &  &  &  &  &  &  & 1 &  &  &  &  &  &  &  \\
 &  &  &  &  &  &  &  &  & 1 &  &  &  &  &  &  \\
 &  &  &  &  &  &  &  &  &  & 1 &  &  &  &  &  \\
 &  &  &  &  &  &  &  &  &  &  & 1 &  &  &  &
\end{array}\right)
$$

$$
\rho((e \otimes e  \otimes e  \otimes j)) =
\left(\begin{array}{rrrr|rrrr|rrrr|rrrr}
 &  &  &  &  &  &  &  & 1 &  &  &  &  &  &  &  \\
 &  &  &  &  &  &  &  &  & 1 &  &  &  &  &  &  \\
 &  &  &  &  &  &  &  &  &  & 1 &  &  &  &  &  \\
 &  &  &  &  &  &  &  &  &  &  & 1 &  &  &  &  \\
\hline
  &  &  &  &  &  &  &  &  &  &  &  & 1 &  &  &  \\
 &  &  &  &  &  &  &  &  &  &  &  &  & 1 &  &  \\
 &  &  &  &  &  &  &  &  &  &  &  &  &  & 1 &  \\
 &  &  &  &  &  &  &  &  &  &  &  &  &  &  & 1 \\
\hline
 -1 &  &  &  &  &  &  &  &  &  &  &  &  &  &  &  \\
 & -1 &  &  &  &  &  &  &  &  &  &  &  &  &  &  \\
 &  & -1 &  &  &  &  &  &  &  &  &  &  &  &  &  \\
 &  &  & -1 &  &  &  &  &  &  &  &  &  &  &  &  \\
\hline
  &  &  &  & -1 &  &  &  &  &  &  &  &  &  &  &  \\
 &  &  &  &  & -1 &  &  &  &  &  &  &  &  &  &  \\
 &  &  &  &  &  & -1 &  &  &  &  &  &  &  &  &  \\
 &  &  &  &  &  &  & -1 &  &  &  &  &  &  &  &
\end{array}\right)
$$
******************************************************************************************\\
it is easily confirmed that the group generated by the following is a solvable group of order $2 \cdot 4^4$,
$$
\begin{aligned}
(i \otimes e \otimes e \otimes e),(j \otimes e \otimes e \otimes e)\\
(e \otimes i \otimes e \otimes e),(e \otimes j \otimes e \otimes e)\\
(e \otimes e \otimes i \otimes e),(e \otimes e \otimes j \otimes e)\\
(e \otimes e \otimes e \otimes i),(e \otimes e \otimes e \otimes j)\\
\end{aligned}
$$
and this group has $257$ conjugacy classes and only one $16$-dimensional irreducible representation.In fact, it is easy to see that the above representation is the irreducible representation.\\\\
In addition, as before, we can construct a $16 \times 16 = 256$-dimensional lattice with the group of ranks $2^{21} \cdot 3^5 \cdot 5^2 \cdot 7$ mentioned above.

\section{Appendix.1}
Using the representation, $\rho(x_i)$ will be
$$
\rho(x_1) = \frac{1}{2}\left(\begin{array}{rrrrrrrrrrrrrrrr}
1 & 0 & 0 & -1 & 0 & 0 & 0 & 0 & 0 & 0 & 0 & 0 & 0 & 1 & 1 & 0 \\
0 & 1 & -1 & 0 & 0 & 0 & 0 & 0 & 0 & 0 & 0 & 0 & 1 & 0 & 0 & 1 \\
0 & -1 & 1 & 0 & 0 & 0 & 0 & 0 & 0 & 0 & 0 & 0 & 1 & 0 & 0 & 1 \\
-1 & 0 & 0 & 1 & 0 & 0 & 0 & 0 & 0 & 0 & 0 & 0 & 0 & 1 & 1 & 0 \\
0 & 0 & 0 & 0 & 1 & 0 & 0 & 1 & 0 & -1 & 1 & 0 & 0 & 0 & 0 & 0 \\
0 & 0 & 0 & 0 & 0 & 1 & 1 & 0 & -1 & 0 & 0 & 1 & 0 & 0 & 0 & 0 \\
0 & 0 & 0 & 0 & 0 & 1 & 1 & 0 & 1 & 0 & 0 & -1 & 0 & 0 & 0 & 0 \\
0 & 0 & 0 & 0 & 1 & 0 & 0 & 1 & 0 & 1 & -1 & 0 & 0 & 0 & 0 & 0 \\
0 & 0 & 0 & 0 & 0 & -1 & 1 & 0 & 1 & 0 & 0 & 1 & 0 & 0 & 0 & 0 \\
0 & 0 & 0 & 0 & -1 & 0 & 0 & 1 & 0 & 1 & 1 & 0 & 0 & 0 & 0 & 0 \\
0 & 0 & 0 & 0 & 1 & 0 & 0 & -1 & 0 & 1 & 1 & 0 & 0 & 0 & 0 & 0 \\
0 & 0 & 0 & 0 & 0 & 1 & -1 & 0 & 1 & 0 & 0 & 1 & 0 & 0 & 0 & 0 \\
0 & 1 & 1 & 0 & 0 & 0 & 0 & 0 & 0 & 0 & 0 & 0 & 1 & 0 & 0 & -1 \\
1 & 0 & 0 & 1 & 0 & 0 & 0 & 0 & 0 & 0 & 0 & 0 & 0 & 1 & -1 & 0 \\
1 & 0 & 0 & 1 & 0 & 0 & 0 & 0 & 0 & 0 & 0 & 0 & 0 & -1 & 1 & 0 \\
0 & 1 & 1 & 0 & 0 & 0 & 0 & 0 & 0 & 0 & 0 & 0 & -1 & 0 & 0 & 1
\end{array}\right)
$$

$$
\rho(x_2) = \frac{1}{2}
\left(\begin{array}{rrrrrrrrrrrrrrrr}
1 & 1 & 1 & -1 & 0 & 0 & 0 & 0 & 0 & 0 & 0 & 0 & 0 & 0 & 0 & 0 \\
1 & 1 & -1 & 1 & 0 & 0 & 0 & 0 & 0 & 0 & 0 & 0 & 0 & 0 & 0 & 0 \\
1 & -1 & 1 & 1 & 0 & 0 & 0 & 0 & 0 & 0 & 0 & 0 & 0 & 0 & 0 & 0 \\
-1 & 1 & 1 & 1 & 0 & 0 & 0 & 0 & 0 & 0 & 0 & 0 & 0 & 0 & 0 & 0 \\
0 & 0 & 0 & 0 & 1 & -1 & 1 & 1 & 0 & 0 & 0 & 0 & 0 & 0 & 0 & 0 \\
0 & 0 & 0 & 0 & -1 & 1 & 1 & 1 & 0 & 0 & 0 & 0 & 0 & 0 & 0 & 0 \\
0 & 0 & 0 & 0 & 1 & 1 & 1 & -1 & 0 & 0 & 0 & 0 & 0 & 0 & 0 & 0 \\
0 & 0 & 0 & 0 & 1 & 1 & -1 & 1 & 0 & 0 & 0 & 0 & 0 & 0 & 0 & 0 \\
0 & 0 & 0 & 0 & 0 & 0 & 0 & 0 & 1 & 1 & -1 & 1 & 0 & 0 & 0 & 0 \\
0 & 0 & 0 & 0 & 0 & 0 & 0 & 0 & 1 & 1 & 1 & -1 & 0 & 0 & 0 & 0 \\
0 & 0 & 0 & 0 & 0 & 0 & 0 & 0 & -1 & 1 & 1 & 1 & 0 & 0 & 0 & 0 \\
0 & 0 & 0 & 0 & 0 & 0 & 0 & 0 & 1 & -1 & 1 & 1 & 0 & 0 & 0 & 0 \\
0 & 0 & 0 & 0 & 0 & 0 & 0 & 0 & 0 & 0 & 0 & 0 & 1 & -1 & -1 & -1 \\
0 & 0 & 0 & 0 & 0 & 0 & 0 & 0 & 0 & 0 & 0 & 0 & -1 & 1 & -1 & -1 \\
0 & 0 & 0 & 0 & 0 & 0 & 0 & 0 & 0 & 0 & 0 & 0 & -1 & -1 & 1 & -1 \\
0 & 0 & 0 & 0 & 0 & 0 & 0 & 0 & 0 & 0 & 0 & 0 & -1 & -1 & -1 & 1
\end{array}\right)
$$

$$
\rho(x_3) = \frac{1}{2}
\left(\begin{array}{rrrrrrrrrrrrrrrr}
1 & 0 & -1 & 0 & 0 & 0 & 0 & 0 & 0 & 0 & 0 & 0 & 0 & 1 & 0 & 1 \\
0 & 1 & 0 & -1 & 0 & 0 & 0 & 0 & 0 & 0 & 0 & 0 & 1 & 0 & 1 & 0 \\
-1 & 0 & 1 & 0 & 0 & 0 & 0 & 0 & 0 & 0 & 0 & 0 & 0 & 1 & 0 & 1 \\
0 & -1 & 0 & 1 & 0 & 0 & 0 & 0 & 0 & 0 & 0 & 0 & 1 & 0 & 1 & 0 \\
0 & 0 & 0 & 0 & 1 & 0 & 1 & 0 & 0 & -1 & 0 & 1 & 0 & 0 & 0 & 0 \\
0 & 0 & 0 & 0 & 0 & 1 & 0 & 1 & -1 & 0 & 1 & 0 & 0 & 0 & 0 & 0 \\
0 & 0 & 0 & 0 & 1 & 0 & 1 & 0 & 0 & 1 & 0 & -1 & 0 & 0 & 0 & 0 \\
0 & 0 & 0 & 0 & 0 & 1 & 0 & 1 & 1 & 0 & -1 & 0 & 0 & 0 & 0 & 0 \\
0 & 0 & 0 & 0 & 0 & -1 & 0 & 1 & 1 & 0 & 1 & 0 & 0 & 0 & 0 & 0 \\
0 & 0 & 0 & 0 & -1 & 0 & 1 & 0 & 0 & 1 & 0 & 1 & 0 & 0 & 0 & 0 \\
0 & 0 & 0 & 0 & 0 & 1 & 0 & -1 & 1 & 0 & 1 & 0 & 0 & 0 & 0 & 0 \\
0 & 0 & 0 & 0 & 1 & 0 & -1 & 0 & 0 & 1 & 0 & 1 & 0 & 0 & 0 & 0 \\
0 & 1 & 0 & 1 & 0 & 0 & 0 & 0 & 0 & 0 & 0 & 0 & 1 & 0 & -1 & 0 \\
1 & 0 & 1 & 0 & 0 & 0 & 0 & 0 & 0 & 0 & 0 & 0 & 0 & 1 & 0 & -1 \\
0 & 1 & 0 & 1 & 0 & 0 & 0 & 0 & 0 & 0 & 0 & 0 & -1 & 0 & 1 & 0 \\
1 & 0 & 1 & 0 & 0 & 0 & 0 & 0 & 0 & 0 & 0 & 0 & 0 & -1 & 0 & 1
\end{array}\right)
$$
$$
\rho(x_4) =
\left(\begin{array}{rrrrrrrrrrrrrrrr}
0 & 0 & 0 & 0 & 0 & 0 & -1 & 0 & 0 & 0 & 0 & 0 & 0 & 0 & 0 & 0 \\
0 & 1 & 0 & 0 & 0 & 0 & 0 & 0 & 0 & 0 & 0 & 0 & 0 & 0 & 0 & 0 \\
0 & 0 & 0 & 0 & -1 & 0 & 0 & 0 & 0 & 0 & 0 & 0 & 0 & 0 & 0 & 0 \\
0 & 0 & 0 & 1 & 0 & 0 & 0 & 0 & 0 & 0 & 0 & 0 & 0 & 0 & 0 & 0 \\
0 & 0 & -1 & 0 & 0 & 0 & 0 & 0 & 0 & 0 & 0 & 0 & 0 & 0 & 0 & 0 \\
0 & 0 & 0 & 0 & 0 & 1 & 0 & 0 & 0 & 0 & 0 & 0 & 0 & 0 & 0 & 0 \\
-1 & 0 & 0 & 0 & 0 & 0 & 0 & 0 & 0 & 0 & 0 & 0 & 0 & 0 & 0 & 0 \\
0 & 0 & 0 & 0 & 0 & 0 & 0 & 1 & 0 & 0 & 0 & 0 & 0 & 0 & 0 & 0 \\
0 & 0 & 0 & 0 & 0 & 0 & 0 & 0 & 1 & 0 & 0 & 0 & 0 & 0 & 0 & 0 \\
0 & 0 & 0 & 0 & 0 & 0 & 0 & 0 & 0 & 0 & 0 & 0 & 0 & 0 & 0 & -1 \\
0 & 0 & 0 & 0 & 0 & 0 & 0 & 0 & 0 & 0 & 1 & 0 & 0 & 0 & 0 & 0 \\
0 & 0 & 0 & 0 & 0 & 0 & 0 & 0 & 0 & 0 & 0 & 0 & 0 & -1 & 0 & 0 \\
0 & 0 & 0 & 0 & 0 & 0 & 0 & 0 & 0 & 0 & 0 & 0 & 1 & 0 & 0 & 0 \\
0 & 0 & 0 & 0 & 0 & 0 & 0 & 0 & 0 & 0 & 0 & -1 & 0 & 0 & 0 & 0 \\
0 & 0 & 0 & 0 & 0 & 0 & 0 & 0 & 0 & 0 & 0 & 0 & 0 & 0 & 1 & 0 \\
0 & 0 & 0 & 0 & 0 & 0 & 0 & 0 & 0 & -1 & 0 & 0 & 0 & 0 & 0 & 0
\end{array}\right)
$$

$$
\rho(x_5) =
\left(\begin{array}{rrrrrrrrrrrrrrrr}
0 & 0 & 0 & 0 & 1 & 0 & 0 & 0 & 0 & 0 & 0 & 0 & 0 & 0 & 0 & 0 \\
0 & 1 & 0 & 0 & 0 & 0 & 0 & 0 & 0 & 0 & 0 & 0 & 0 & 0 & 0 & 0 \\
0 & 0 & 0 & 0 & 0 & 0 & -1 & 0 & 0 & 0 & 0 & 0 & 0 & 0 & 0 & 0 \\
0 & 0 & 0 & 1 & 0 & 0 & 0 & 0 & 0 & 0 & 0 & 0 & 0 & 0 & 0 & 0 \\
1 & 0 & 0 & 0 & 0 & 0 & 0 & 0 & 0 & 0 & 0 & 0 & 0 & 0 & 0 & 0 \\
0 & 0 & 0 & 0 & 0 & 1 & 0 & 0 & 0 & 0 & 0 & 0 & 0 & 0 & 0 & 0 \\
0 & 0 & -1 & 0 & 0 & 0 & 0 & 0 & 0 & 0 & 0 & 0 & 0 & 0 & 0 & 0 \\
0 & 0 & 0 & 0 & 0 & 0 & 0 & 1 & 0 & 0 & 0 & 0 & 0 & 0 & 0 & 0 \\
0 & 0 & 0 & 0 & 0 & 0 & 0 & 0 & 1 & 0 & 0 & 0 & 0 & 0 & 0 & 0 \\
0 & 0 & 0 & 0 & 0 & 0 & 0 & 0 & 0 & 0 & 0 & 0 & 0 & 1 & 0 & 0 \\
0 & 0 & 0 & 0 & 0 & 0 & 0 & 0 & 0 & 0 & 1 & 0 & 0 & 0 & 0 & 0 \\
0 & 0 & 0 & 0 & 0 & 0 & 0 & 0 & 0 & 0 & 0 & 0 & 0 & 0 & 0 & -1 \\
0 & 0 & 0 & 0 & 0 & 0 & 0 & 0 & 0 & 0 & 0 & 0 & 1 & 0 & 0 & 0 \\
0 & 0 & 0 & 0 & 0 & 0 & 0 & 0 & 0 & 1 & 0 & 0 & 0 & 0 & 0 & 0 \\
0 & 0 & 0 & 0 & 0 & 0 & 0 & 0 & 0 & 0 & 0 & 0 & 0 & 0 & 1 & 0 \\
0 & 0 & 0 & 0 & 0 & 0 & 0 & 0 & 0 & 0 & 0 & -1 & 0 & 0 & 0 & 0
\end{array}\right)
$$

$$
\rho(x_6) =\frac{1}{2}
\left(\begin{array}{rrrrrrrrrrrrrrrr}
1 & 0 & 0 & 0 & 0 & 1 & 0 & 0 & 0 & 0 & 0 & -1 & 0 & 0 & -1 & 0 \\
0 & 1 & 0 & 0 & -1 & 0 & 0 & 0 & 0 & 0 & -1 & 0 & 0 & 0 & 0 & 1 \\
0 & 0 & 1 & 0 & 0 & 0 & 0 & 1 & 0 & 1 & 0 & 0 & 1 & 0 & 0 & 0 \\
0 & 0 & 0 & 1 & 0 & 0 & -1 & 0 & 1 & 0 & 0 & 0 & 0 & -1 & 0 & 0 \\
0 & -1 & 0 & 0 & 1 & 0 & 0 & 0 & 0 & 0 & -1 & 0 & 0 & 0 & 0 & 1 \\
1 & 0 & 0 & 0 & 0 & 1 & 0 & 0 & 0 & 0 & 0 & 1 & 0 & 0 & 1 & 0 \\
0 & 0 & 0 & -1 & 0 & 0 & 1 & 0 & 1 & 0 & 0 & 0 & 0 & -1 & 0 & 0 \\
0 & 0 & 1 & 0 & 0 & 0 & 0 & 1 & 0 & -1 & 0 & 0 & -1 & 0 & 0 & 0 \\
0 & 0 & 0 & 1 & 0 & 0 & 1 & 0 & 1 & 0 & 0 & 0 & 0 & 1 & 0 & 0 \\
0 & 0 & 1 & 0 & 0 & 0 & 0 & -1 & 0 & 1 & 0 & 0 & -1 & 0 & 0 & 0 \\
0 & -1 & 0 & 0 & -1 & 0 & 0 & 0 & 0 & 0 & 1 & 0 & 0 & 0 & 0 & 1 \\
-1 & 0 & 0 & 0 & 0 & 1 & 0 & 0 & 0 & 0 & 0 & 1 & 0 & 0 & -1 & 0 \\
0 & 0 & 1 & 0 & 0 & 0 & 0 & -1 & 0 & -1 & 0 & 0 & 1 & 0 & 0 & 0 \\
0 & 0 & 0 & -1 & 0 & 0 & -1 & 0 & 1 & 0 & 0 & 0 & 0 & 1 & 0 & 0 \\
-1 & 0 & 0 & 0 & 0 & 1 & 0 & 0 & 0 & 0 & 0 & -1 & 0 & 0 & 1 & 0 \\
0 & 1 & 0 & 0 & 1 & 0 & 0 & 0 & 0 & 0 & 1 & 0 & 0 & 0 & 0 & 1
\end{array}\right)
$$

$$
\rho(x_7) =\frac{1}{2}
\left(\begin{array}{rrrrrrrrrrrrrrrr}
1 & 0 & 0 & 0 & 0 & 1 & 0 & 0 & 0 & 1 & 0 & 0 & 1 & 0 & 0 & 0 \\
0 & 1 & 0 & 0 & 1 & 0 & 0 & 0 & 1 & 0 & 0 & 0 & 0 & 1 & 0 & 0 \\
0 & 0 & 1 & 0 & 0 & 0 & 0 & 1 & 0 & 0 & 0 & -1 & 0 & 0 & -1 & 0 \\
0 & 0 & 0 & 1 & 0 & 0 & 1 & 0 & 0 & 0 & -1 & 0 & 0 & 0 & 0 & -1 \\
0 & -1 & 0 & 0 & 1 & 0 & 0 & 0 & 1 & 0 & 0 & 0 & 0 & -1 & 0 & 0 \\
-1 & 0 & 0 & 0 & 0 & 1 & 0 & 0 & 0 & 1 & 0 & 0 & -1 & 0 & 0 & 0 \\
0 & 0 & 0 & -1 & 0 & 0 & 1 & 0 & 0 & 0 & -1 & 0 & 0 & 0 & 0 & 1 \\
0 & 0 & -1 & 0 & 0 & 0 & 0 & 1 & 0 & 0 & 0 & -1 & 0 & 0 & 1 & 0 \\
0 & -1 & 0 & 0 & -1 & 0 & 0 & 0 & 1 & 0 & 0 & 0 & 0 & 1 & 0 & 0 \\
-1 & 0 & 0 & 0 & 0 & -1 & 0 & 0 & 0 & 1 & 0 & 0 & 1 & 0 & 0 & 0 \\
0 & 0 & 0 & 1 & 0 & 0 & 1 & 0 & 0 & 0 & 1 & 0 & 0 & 0 & 0 & 1 \\
0 & 0 & 1 & 0 & 0 & 0 & 0 & 1 & 0 & 0 & 0 & 1 & 0 & 0 & 1 & 0 \\
-1 & 0 & 0 & 0 & 0 & 1 & 0 & 0 & 0 & -1 & 0 & 0 & 1 & 0 & 0 & 0 \\
0 & -1 & 0 & 0 & 1 & 0 & 0 & 0 & -1 & 0 & 0 & 0 & 0 & 1 & 0 & 0 \\
0 & 0 & 1 & 0 & 0 & 0 & 0 & -1 & 0 & 0 & 0 & -1 & 0 & 0 & 1 & 0 \\
0 & 0 & 0 & 1 & 0 & 0 & -1 & 0 & 0 & 0 & -1 & 0 & 0 & 0 & 0 & 1
\end{array}\right)
$$

\section{Appendix.2}
The irreducible representation of the quaternion is as follows.
$$
\begin{aligned}
&\boldsymbol{i} := \tau(i) =\left(\begin{array}{cccc}
 & 1 &  &  \\
-1 &  &  &  \\
 &  &  & 1 \\
 &  & -1 &
\end{array}\right) \\\\
&\boldsymbol{j} :=\tau(j) =\left(\begin{array}{cccc}
 &  & 1 &  \\
 &  &  & -1 \\
-1 &  &  &  \\
 & 1 &  &
\end{array}\right) \\\\
&\boldsymbol{k}:=\tau(k) =\left(\begin{array}{cccc}
 &  &  & -1 \\
 &  & -1 &  \\
 & 1 &  &  \\
1 &  &  &
\end{array}\right)
\end{aligned}
$$
This allows the alternation group $2 \cdot A_4$ described above to be written as follows.
$$
\left\langle\left(\begin{array}{rrrr}
 & 1 &  &  \\
-1 &  &  &  \\
 &  &  & 1 \\
 &  & -1 &
\end{array}\right),\left(\begin{array}{rrrr}
 &  & 1 &  \\
 &  &  & -1 \\
-1 &  &  &  \\
 & 1 &  &
\end{array}\right), \frac{1}{2} \left(\begin{array}{rrrr}
-1 & 1 & 1 & -1 \\
-1 & -1 & -1 & -1 \\
-1 & 1 & -1 & 1 \\
1 & 1 & -1 & -1
\end{array}\right)\right\rangle \cong 2 \cdot A_4
$$

\end{document}